\documentclass[11pt,a4paper]{article}
\usepackage[T1]{fontenc}
\usepackage[utf8]{inputenc}
\usepackage{lmodern,amsmath,amssymb,amsthm,microtype,mathtools}
\usepackage[a4paper,left=26mm,right=26mm,top=23mm,bottom=24mm,
            headheight=14pt,headsep=7mm,footskip=11mm]{geometry}
\usepackage{enumitem,caption,fancyhdr,needspace,booktabs,amssymb}
\usepackage[section]{placeins}
\usepackage{tikz}
\usetikzlibrary{arrows.meta,calc}
\usepackage[unicode,bookmarksnumbered]{hyperref}
\usepackage[capitalise, noabbrev]{cleveref}
\definecolor{inkblue}{HTML}{214B68}
\definecolor{routeblue}{HTML}{246489}
\definecolor{partgreen}{HTML}{26695D}
\definecolor{warm}{HTML}{98502D}
\definecolor{lightline}{HTML}{BDCAD0}
\hypersetup{colorlinks,linkcolor=inkblue,citecolor=inkblue,urlcolor=inkblue,
  pdftitle={The Crossing Number of Crossing-Critical Graphs Revisited},
  pdfauthor={Ilkyoo Choi, John Baptist Gauci, Petr Hliněný, Tomáš Masařík, Alexandra Wesolek}}
\newtheorem{theorem}{Theorem}
\newtheorem{lemma}[theorem]{Lemma}
\newtheorem{corollary}[theorem]{Corollary}
\newcommand{\crn}{\operatorname{cr}}
\newcommand{\Def}{\operatorname{def}}
\newcommand{\cI}{\mathcal I}
\newcommand{\cX}{\mathcal X}
\newcommand{\step}[1]{\par\smallskip\noindent\textit{#1}\enspace\ignorespaces}
\newcommand{\doi}[1]{\href{https://doi.org/#1}{doi:\,#1}}
\setlist[enumerate]{label=\textup{(\roman*)},leftmargin=2.2em,itemsep=2pt,topsep=4pt}
\tikzset{vertex/.style={circle,fill=black,inner sep=1.8pt},
  block/.style={circle,draw=inkblue,fill=white,line width=.85pt,inner sep=2.15pt},
  edge/.style={line width=.85pt,line cap=round,line join=round},
  lab/.style={font=\small}, panel/.style={font=\footnotesize\sffamily,text=inkblue}}
\fancypagestyle{first}{\fancyhf{}\fancyfoot[C]{\small\thepage}}
\begin{document}
\thispagestyle{first}
\begin{center}
  {\LARGE\bfseries The Crossing Number of Crossing-Critical Graphs Revisited\par}
  \vspace{9pt}
  {\large Ilkyoo Choi\footnote{Department of Mathematics, Hankuk University of Foreign Studies, Republic of Korea.
  Discrete Mathematics Group, Institute for Basic Science (IBS), Republic of Korea.  
  Center for AI and Natural Sciences, Korea Institute for Advanced Study (KIAS), Republic of Korea.},\quad John Baptist Gauci\footnote{Department of Mathematics, Faculty of Science, University of Malta, Malta},\quad Petr Hlin\v{e}n\'y\footnote{Faculty of Informatics, Masaryk University, Brno, Czech Republic},\\[3pt]
    Tom\'a\v{s} Masa\v{r}\'ik\footnote{University of Warsaw, Poland},\quad Alexandra Wesolek\footnote{University of Bordeaux, CNRS, LaBRI, Bordeaux, France}\par}
\end{center}

\begin{abstract}
We prove that for any positive integer $k$, every \(k\)-crossing-critical graph has crossing number at most \(\tfrac32 k+c\) for an absolute constant \(c\) independent of $k$.
This is an improvement over the original upper bound of $\frac52k+16$ by Richter and Thomassen, as well as the fairly recently improved bound of $2k+8\sqrt k+47$ by Bar\'at and T\'oth.
Our bound is the first one which breaks the barrier of $2k+o(k)$ estimates for general graphs.
\end{abstract}

\section{Introduction}\label{sec:intro}

The crossing number \(\crn(G)\) is the least number of crossings in a drawing of \(G\) in the plane. We refer to \Cref{sec:min} for technical details.

A graph \(G\) is \emph{\(k\)-crossing-critical} if \(\crn(G)\ge k\) and \(\crn(G-e)<k\) for every edge~\(e\). 
While the $1$-crossing-critical graphs are exactly the subdivisions of graphs $K_5$ and $K_{3,3}$, already the early constructions of \v{S}ir\'a\v{n} \cite{siran84} and Kochol \cite{kochol87} show that the structure of $k$-crossing-critical graphs is quite rich for any $k\ge2$.
Nowadays, the $2$-crossing-critical graphs are almost completely described in the work of Bokal, Oporowski, Richter and Salazar \cite{bokalOporowskiRichter16}, and an implicit characterisation of the classes of $k$-crossing-critical graphs for very $k\ge3$ was given by Dvo\v{r}\'ak, Hlin\v{e}n\'y and Mohar \cite{dvorakHlinenyMohar18}.

Another early influential result on crossing-critical graphs was by Richter and Thomassen \cite{DBLP:journals/jct/RichterT93}, who proved that $k$-crossing-critical graphs have their crossing number bounded in terms of $k$.
In this respect, note that the definition of criticality does not a priori bound $\crn(G)$ from above, and it is not immediately clear whether this should be bounded at all since constructions show possibility of an arbitrarily large additive gap between $k$ and $\crn(G)$, of order up to $\Theta(\sqrt k)$ (consider, e.g., the complete or complete bipartite graphs).

How far above \(k\) can \(\crn(G)\) be? The aforementioned result \cite{DBLP:journals/jct/RichterT93} by Richter and Thomassen shows $\crn(G)\leq2.5k+16$. Balogh, Lea\~nos and Salazar \cite{BLS} later obtained a bound of the form \(\tfrac32 k+O(1)\) for the special case of large sparse graphs in which the number of edges is not artificially inflated by operations such as edge subdivisions.
Bar\'at and T\'oth~\cite{BT} improved the general upper bound to $2k+8\sqrt k+47$, and more recently Ding, Ouyang, Huang, and Dong \cite{DBLP:journals/dm/DingOHD22} proved further slight improvements of the latter bound depending on minimum degree of the graph.

The aforementioned characterisations of $k$-crossing-critical graphs \cite{bokalOporowskiRichter16,dvorakHlinenyMohar18} do not seem to help with answering this question since they do not provide a description of the finite sets of ``basic graphs'' for each $k$.

We obtain the following general bound, essentially matching the special bound of Balogh, Lea\~nos and Salazar \cite{BLS}, but for all crossing-critical graphs without further assumptions.

\begin{theorem}\label{thm:main}
There is an absolute constant \(c\) such that, for every integer \(k\ge1\), if a simple graph \(G\) is $k$-crossing critical, then $\crn(G)\le\tfrac32 k+c$.
\end{theorem}

We prove the theorem for finite loopless multigraphs. Parallel edges arise naturally when a vertex
with two neighbours is suppressed. Degrees and edge counts in \(G\) include multiplicities. The
\emph{simplification} \(\widetilde{G}\) is obtained from \(G\) by keeping one edge from each parallel
class. Thus \(d_{\widetilde{G}}(v)\) is the number of distinct neighbours of \(v\) in \(G\).

There is also a deletion form of \Cref{thm:main}. If \(G\) has at least one edge, take
\(k=1+\max_e\crn(G-e)\). An edge attaining the maximum then satisfies $\crn(G-e)\ge\tfrac23\bigl(\crn(G)-c\bigr)-1$.

\paragraph{Idea of the proof.}
The factor $3/2$ comes from special edges $e=\{u,v\}\in E(G)$, which have the special property that their endpoints $u,v$ in $G-e$ are connected by two short edge-disjoint paths and the internal vertices on the paths have small degree. If such a special edge exists in the graph, we can redraw $e$ along one of the short paths in an optimal drawing of $G-e$. If $G-e$ has at most $k-1$ crossings, one of the paths is crossed at most $(k-1)/2$ times, modulo technical details. Since the path lengths and the degrees of their internal vertices are bounded, tracing that path gives a drawing of \(G\) with at most
\(3(k-1)/2+O(1)\) crossings.

In order to show that such special edges must exist, we first show that the number of edges and crossings is at most linear in $n=V(G)$ (\Cref{lem:min}). We use the planar separator theorem to remove a small separator $B$ with at most $\varepsilon n$ vertices, to which we also add all high degree vertices, such that the vertices in $G-B$ have bounded degree in $G$ and form small components (\Cref{lem:separator}). If a component contains a cycle with a chord, we have found a special edge (the chord). If not, we use a degree bound for chordless graphs which forces many edges between $B$ and $U$. Using these edges, we essentially find a different configuration (a subdivided special edge), which allows us to delete and reinsert an edge at a low cost.

\section{A minimum counterexample}\label{sec:min}

In a \emph{drawing} of a graph~$G$ in the plane, vertices are mapped to distinct points and each edge is represented by a simple continuous arc connecting its corresponding endpoints. We require that edges do not pass through other vertices and that no three edges meet in a common interior point. Further, we assume that edges do not touch and that they intersect in at most a finite number of points. A \emph{crossing} in a drawing~$D$ is a common interior point of two edges of~$D$. A drawing~$D$ is \emph{simple} if adjacent edges do not cross and non-adjacent edges cross at most once. Every graph has a simple optimal drawing that exhibits the least number of crossings~\cite{PT}.

Parallel edges may be drawn together in a \emph{bundle}. To see this, start from a maximum subset of parallel edges in a simple optimal drawing. One of the parallel edges is crossed the least number of times, say $t$ times, and we can redraw the other edges such that they are also only crossed $t$ times each. 

From the following lemma we can deduce an easy upper bound on the number of crossings of $k$-critical graphs in terms of $k$ and their number of edges. 

\begin{lemma}\label{lem:insertion}
For each edge $e$ of a loopless multigraph $G$,
$\crn(G)\le \crn(G-e)+|E(G)|-1$.
\end{lemma}
The lemma can be obtained by deleting $e$, taking an optimal drawing of $G-e$ and redrawing by adding $e$, see e.g. Pach and T\'oth~\cite{PT}. This simple lemma can be used to reduce to the case where the number of edges and crossings in $G$ are at most linear. 

\begin{lemma}\label{lem:min}
Fix $\alpha>1$ and an integer \(k\ge1\). There exist constants $c_1,c_2>0$, depending only on $\alpha$, such that the following is true.
Among the loopless multigraphs satisfying $\crn(G)>\alpha k$ and $\crn(G-e)\le k-1$ for every $e\in E(G)$, choose one of minimum order, if any exists. 
Then $\delta(\widetilde{G})\ge3$, $|E(G)|\le c_1n$, and $\crn(G)\le c_2n$.

In particular, one can choose $c_1=40$ and $c_2=120$ for $\alpha=\frac32$.
\end{lemma}
\begin{proof}
Write \(n=|V(G)|\), \(m=|E(G)|\) and \(x=\crn(G)\). Note that $x>\alpha k>0$, so $m>0$ as well.

\step{Three distinct neighbours.}
A vertex with at most one neighbour can be removed and reinserted without adding a crossing. If it
has exactly one neighbour, the vertex and  all its incident edges can be  redrawn in a drawing without $v$, which has at least one edge less and hence we can draw it with at most $k-1$ crossings. Hence, minimality excludes such a vertex.

Suppose \(z\) has just two neighbours \(u,v\), with \(a\) edges to \(u\) and \(b\) to \(v\), where
\(a\le b\). Replace \(z\) and its incident edges by \(a\) new \(uv\)-edges. Call the result \(G'\) (Figure~\ref{fig:suppress}).
In a bundled optimal drawing of \(G'\), subdivide one of the $a$ new edges close to $v$ and reroute the other $a-1$ edges so their new endpoint is the vertex $z'$ obtained from subdivision. Now draw $b-1$ edges close to the $z'v$ edge in a bundle. The final drawing has the same number of edges as the drawing we started with. This gives \(\crn(G')\ge\crn(G)\).

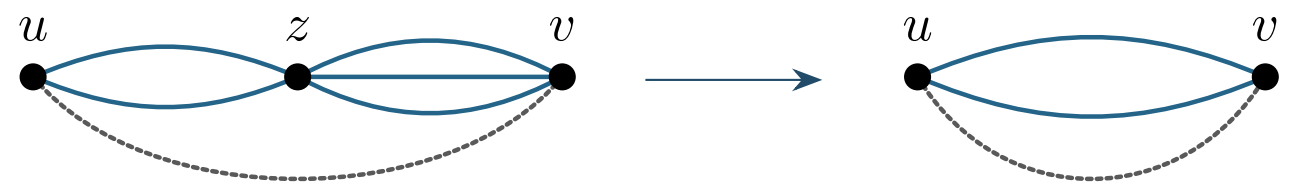
\begin{figure}[!htbp]
\centering
\begin{tikzpicture}[x=1cm,y=1cm,lab/.append style={font=\small}]
  \coordinate (U) at (0,0); \coordinate (Z) at (1.75,0); \coordinate (V) at (3.5,0);
  \draw[edge,black!65,densely dotted] (U) .. controls (.75,-.90) and (2.75,-.90) .. (V);
  \draw[edge,routeblue] (U) to[bend left=23] (Z) (U) to[bend right=23] (Z);
  \draw[edge,routeblue] (Z) to[bend left=28] (V) (Z)--(V) (Z) to[bend right=28] (V);
  \foreach \point/\letter in {U/u,Z/z,V/v}{
    \node[vertex] at (\point) {}; \node[lab,above=3pt] at (\point) {$\letter$};}
  \draw[-{Stealth[length=2mm]},inkblue] (4.05,-.02)--(5.22,-.02);
  \coordinate (U1) at (5.85,0); \coordinate (V1) at (8.15,0);
  \draw[edge,black!65,densely dotted] (U1) .. controls (6.4,-.90) and (7.6,-.90) .. (V1);
  \draw[edge,routeblue] (U1) to[bend left=23] (V1) (U1) to[bend right=23] (V1);
  \node[vertex] at (U1) {}; \node[vertex] at (V1) {};
  \node[lab,above=3pt] at (U1) {$u$}; \node[lab,above=3pt] at (V1) {$v$};
\end{tikzpicture}
\caption{Suppressing \(z\) when the multiplicities are two and three. The two new edges follow the
smaller bundle. The dotted \(uv\)-edge was already present.}
\label{fig:suppress}
\end{figure}

To delete one of the new edges, first delete one \(uz\)-edge in \(G\) and apply the replacement with multiplicities \(a-1,b\) for the optimal drawing of the resulting graph. 
If \(a=1\), simply remove the resulting vertex. Thus \(G'\) satisfies every
deletion inequality for the same \(k\). Its smaller order contradicts minimality.

\step{Parallel classes have size at most two.}
From a class of size \(a\ge3\), delete one edge \(e\). In a simple optimal drawing of \(G-e\), the
remaining \(a-1\) copies never cross one another. Together they account for at most \(k-1\)
crossings. Drawing \(e\) beside the least-crossed copy gives
$x\le (k-1)+(k-1)/(a-1)\le\tfrac32(k-1)$,
a contradiction. Thus \(s\le m\le2s\), where \(s=|E(\widetilde{G})|\).

\step{Enumerating edges and crossings.}
Since there is an edge in $G$, Lemma~\ref{lem:insertion} implies \(x<k+m\). 
Also \(x>\alpha k\), so $(\alpha-1)k<m$ and $x<k+m<\frac{\alpha}{\alpha-1}m$. Altogether, \(x<\frac{\alpha}{\alpha-1}m\le\frac{2\alpha}{\alpha-1}s\).

Finally, we apply a short probabilistic argument.
In a simple optimal drawing of \(\widetilde{G}\), we keep each vertex independently with probability \(p=1/4\) and take the induced drawing.
A vertex, an edge and a crossing survive with probabilities \(p,p^2,p^4\), respectively, since an edge has two ends and a crossing has four distinct ends. 
So, the expected number of surviving crossings is $p^4\crn(\widetilde{G})\le p^4x<\frac{2\alpha}{\alpha-1}p^4s$.
Removing at most one edge for each surviving crossing leaves a simple planar graph $G_1$. 
Euler's bound on $G_1$ yields $|E(G_1)|\leq3|V(G_1)|$ which, after adding back the removed edges, by linearity of expectation gives $p^2s\le3pn+\frac{2\alpha}{\alpha-1}p^4s$.

Substituting \(p=1/4\) in the latter inequality gives \(\frac1{16}s\le\frac34n+\frac{\alpha}{128(\alpha-1)}s\), or \(s\le 96\frac{7\alpha-8}{\alpha-1}n\).
It follows that \(m\le2s\le 192\frac{7\alpha-8}{\alpha-1}n\) and \(x<3m\le 576\frac{7\alpha-8}{\alpha-1}n\).
When $\alpha=\frac32$, we specifically get $m\leq192n/5\leq40n$ and $x<3m\leq 120n$.
\end{proof}

\section{Small components and their neighbours}\label{sec:inc}

To follow our proof path, we need a separator \(B\) in our graph whose deletion leaves components of constant size and whose vertices have constant degree. 
Furthermore, later in \cref{sec:paths}, we will count the neighbours in \(B\) of disjoint connected subsets of those components. 

For a family \(\cX\) of such disjoint sets, we write \(\iota_B(\cX)\) denoting  the number of pairs \((X,b)\) with \(X\in\cX\), \(b\in B\) and at least one edge between \(X\) and \(b\). 
Several edges joining \(X\) to \(b\) thus still contribute only one such incidence. 
A vertex set is called \emph{connected} if it induces a connected subgraph.

\begin{lemma}
\label{lem:separator}
Let $c_1,c_2>1$ be fixed. For each \(0<\varepsilon\le1\), there is an integer \(r_\varepsilon\ge2\) with the following property. Let
\(G\) be a multigraph on \(n\) vertices with at most \(c_1n\) edges and crossing number at
most \(c_2n\). There is a set \(B\subseteq V(G)\) and \(U\coloneqq V(G)\setminus B\) such that: 
\begin{enumerate}
\item \(|B|\le\varepsilon n\),
\item each component of \(G[U]\) has at most \(r_\varepsilon\) vertices, and every \(u\in U\) has \(d_G(u)\le r_\varepsilon\), and
\item every family \(\cX\) of pairwise disjoint non-empty connected sets in \(U\) satisfies
\begin{equation}
\iota_B(\cX)\le 2|\cX|+\varepsilon n . \label{eq:incidence}
\end{equation}
\end{enumerate}
\end{lemma}
\begin{proof}
Put \(\rho\coloneqq\varepsilon/(4c_2)\). For the sampling argument below, we choose \(B\) with \(|B|\le\varepsilon\rho n/4< \varepsilon n\).
\step{Constructing the separator \(B\).}
Planarize a simple optimal drawing of \(\widetilde{G}\) by putting a vertex at each crossing. The resulting
plane graph has \(n_{\mathrm p}=n+\crn(\widetilde{G})\le (c_2+1)n\leq 2c_2n\) vertices. If a component has order \(a>r_\varepsilon\), the planar separator theorem of Lipton and Tarjan~\cite{LT} gives a set of at most \(\lambda\sqrt a\) vertices whose removal leaves components of order at most \(2a/3\), where \(\lambda=2\sqrt{2}\). 
Repeat until every component has order at most \(r_\varepsilon\), and denote by \(Z\) the set of all vertices removed.

We use a standard argument to bound the size of \(Z\). 
When a component of order \(a\) is separated, charge \(\lambda/\sqrt a\) to each of its vertices. A vertex is charged along one chain of nested components. Its last charge is at most \(\lambda/\sqrt{r_\varepsilon}\), and each earlier charge is at most \(\sqrt{2/3}\) times the next. Summing this geometric series gives at most \(\lambda_0/\sqrt{r_\varepsilon}\) per vertex, where \(\lambda_0=\lambda/(1-\sqrt{2/3})= 2\sqrt{2}/(1-\sqrt{2/3}) \). Thus \(|Z|\le \lambda_0 n_{\mathrm p}/\sqrt{r_\varepsilon}\).

We put in \(B\) both the original vertices of \(Z\) as well as the endpoints of both edges at every crossing in
\(Z\), and all vertices of \(G\) of degree greater than \(r_\varepsilon\). The sum of degrees in \(G\) is at
most \(2c_1n\), so $|B|\le4|Z|+\frac{2c_1n}{r_\varepsilon}
\le \left(\frac{8c_2\lambda_0}{\sqrt{r_\varepsilon}}+\frac{2c_1}{r_\varepsilon}\right)n$.
The coefficient at \(n\) tends monotonically to zero as \(r_\varepsilon\) increases. Choose least \(r_\varepsilon\ge2\) so that this coefficient is at most \(\varepsilon\rho/4\). 
Hence, $r_\varepsilon\ge(2\cdot8^2\cdot4)c_1c_2^2\lambda_0^2/(\varepsilon\rho)=2048c_1c_2^3\lambda_0^2/\varepsilon^2$ suffices. Then 
\begin{align}\label{eq:B}
|B|\le\varepsilon\rho n/4<\varepsilon n.
\end{align} 
By definition, an edge with both ends outside \(B\) is represented in the drawing of \(\widetilde{G}\) by a path avoiding \(Z\). 
Therefore, each component of \(G[U]\) lies within a single remaining component of the planarization. This proves (i) and (ii).

\step{Counting the incidences.}
Fix \(\cX=\{X_1,\ldots,X_q\}\) as in (iii). In a simple optimal drawing of~\(G\), planarise each subdrawing of $G[X_i]$, $1\le i\le q$, and keep a spanning tree of this planarisation and exactly one edge from it to each of the neighbours of $X_i$ in \(B\). Denote the resulting graph subdrawing by~$W_i$.
Observe that \(W_i-B\) is connected.
For any $\cI\subseteq\{1,\ldots,q\}$, denote by $W_\cI$ the subdrawing $\bigcup_{i\in\cI}W_i$ (which includes the incident vertices of $B$).
For \(i\neq j\in\cI\), \(W_i\)~does not share vertices with \(W_j\) except within $B$, but their edges can cross in $W_\cI$.

Let \(q_*=|\cI|\), let \(e_*\) count the number of edges of $W_\cI$ incident to the vertices in $B$, and let $x_*$ be the number of crossings in $W_\cI$. We claim that
\begin{align}\label{eq:estar}
e_*\le 2\bigl(|B|+q_*+x_*\bigr).
\end{align}

Make $W_\cI'$ by subdividing edges of $W_\cI$ such that no edge in the subdivision has more than one crossing and no crossed edge has an end vertex from some $W_i$ or $B$ (put this simply, every segment between two crossings in $W_\cI$ or between a vertex and a crossing in $W_\cI$ is subdivided).
Make $W_\cI''$ from $W_\cI'$ by removing any one edge from each crossing, removing at most $x_*$ edges in total. No edge incident to $B$ is removed.
Observe that, as $W_\cI'-B$ has $q_*$ connected components, $W_\cI''-B$ will have at most $q_*+x_*$ connected components.

Since $W_\cI''$ is planar, applying Euler's formula to the simple bipartite planar graph $Q$ obtained from $W_\cI''$ by contracting each of the said components into a vertex gives $|E(Q)|\le 2|V(Q)|$.
Now, as $|E(Q)|=e_*$ and $|V(Q)|\leq|B|+q_*+x_*$, inequality \eqref{eq:estar} follows.

In the remaining argument, we pick every element of $\cI\subseteq\{1,\ldots,q\}$ independently with probability \(\rho\).
We hence have $q_*=|\cI|=\rho q$ in expectation. The total number of edges incident to $B$ over all $W_i$ is by definition equal $\iota_B(\cX)$, and so $e_*=\rho\,\iota_B(\cX)$ in expectation.
And since the crossings in $W_\cI$ come only from crossing pairs of $W_i$ and $W_j$ for \(i\neq j\in\cI\), we also have $x_*\le\rho^2\crn(G)$ in expectation.
By linearity of expectation we can thus from \eqref{eq:estar} conclude $\rho\,\iota_B(\cX)\le 2\bigl(|B|+\rho q+\rho^2\crn(G)\bigr)$, or $\iota_B(\cX)\le 2\bigl(|B|/\rho+q+\rho\crn(G)\bigr)$.
The assumption of the lemma provides $\crn(G)\leq c_2n$, so \(\rho\crn(G)=\varepsilon/(4c_2)\crn(G)\le\varepsilon n/4\).
By \eqref{eq:B}, \(|B|/\rho\le\varepsilon n/4\). Putting things together, $\iota_B(\cX)\le 2(\varepsilon n/4+q+\varepsilon n/4)=\varepsilon n+2q$, which confirms \eqref{eq:incidence} for every choice of $\cX$ and hence (iii) is proven.
\end{proof}

From the separator constructed in Lemma~\ref{lem:separator} and its properties, we can readily improve the asymptotics of the $2k+8\sqrt{k}+47$ bound from~\cite{BT}. The proof of the following corollary also gives an intuition on why the separator lemma is a useful tool in general.
\begin{corollary}
    There is an absolute constant $c$ such that, for every integer $k\geq 1$, if a simple graph $G$ satisfies $\crn(G-e)<k$ for every edge $e$, then
$\crn(G)\le 2 k+c$.
\end{corollary}
\begin{proof}
We choose a minimum-order graph $G$ with $\crn(G)\geq \frac{3}{2}k+c$. 
By Lemma~\ref{lem:min}, we can assume $G$ has minimum degree $3$, at most $40n$ edges and $\crn(G)\leq 120n$.
Fix \(\varepsilon=1/3\). We take $r_\varepsilon$ and $B$ as in Lemma~\ref{lem:separator} and let $U=V(G)\setminus B$. We claim that we can find a cycle $C$ with $|C\cap B|\leq 1$. We then use a standard technique already used in~\cite{DBLP:journals/jct/RichterT93,BT}. We delete an edge $e$ of $C$ such that we get a path $P$ whose inner vertices are in $U$. We then use an optimal drawing of $G-e$ with less than $k$ crossings and add back the edge $e$ tightly along the path $P$. Doing this, we add at most $r_\varepsilon^2$ crossings with edges incident to inner vertices of $P$ and at most $k-1$ crossings with edges that are not on $P$ but cross an edge of $P$. This gives at most $2k+c$ crossings for $c=r_\varepsilon^2-1$.

If there is a cycle in $U$, we are done. Hence we can suppose that each connected component $A$ of $U$ is a tree. Note that the number of edges in $A$ is $|A|-1$. Let $U$ consist of $h$ components. Then the number of edges within $U$ is exactly $|U|-h$. Since $\varepsilon=1/3$ note that $|U|\geq \frac{2}{3}n$. Since every vertex of $U$ has degree at least $3$, the number of edges between $U$ and $B$ is at least
$\sum_{u \in U} \deg(u) - 2(|U|-h) \geq 3|U| - 2(|U|-h) \geq |U|+2h$.
By Lemma~\ref{lem:separator}, if $\cX$ is the family of connected components of $U$, $\iota_B(\cX)\le 2|\cX|+ |B| = 2h+|B|$. Since $|B|+2h < |U|+2   h$, there exists an incidence between some vertex $b\in B$ and some component $A\in\cX$ which accounts for at least $2$ edges between $U$ and $B$. But that means $b$ has two neighbours in $A$ and we find a cycle within $U\cup b$, and we are done by the rerouting argument at the beginning of the proof.
\end{proof}

\section{Reinserting an edge}\label{sec:local}
If we find nice substructures in $G$, such as small cycles whose neighbourhood is small, we can draw $G$ with relatively few crossings based on the drawing of $G-e$ for some edge $e$. This section discusses such redrawing techniques. 
\subsection{Two short paths}\label{sec:twopaths}

\begin{lemma}
\label{lem:traces}
For each integer \(d\ge2\) there is a constant \(a_d\) with the following property. Let \(f\) be
an edge of a loopless multigraph \(G\). Suppose one of the following holds.
\begin{enumerate}
\item \(f=ub\) and there are two edge-disjoint \(u\)--\(b\) paths in \(G-f\).
\item \(f=vb\), where \(v\) is of degree $3$ with distinct neighbours \(b,p,r\), and \(d_G(p)\le d\).
There are two edge-disjoint \(p\)--\(b\) paths in \(G-v\).
\end{enumerate}
If both paths have length at most \(d+2\) and all their internal vertices have degree at most
\(d\) in~\(G\), then
$\crn(G)\le \tfrac32\crn(G-f)+a_d$. 
\end{lemma}
\begin{proof}
Take a simple optimal drawing of \(G-f\) with \(\crn(G-f)\) crossings. 
Denote the two paths from (i) or (ii) by \(P_1,P_2\), and choose $P\in\{P_1,P_2\}$ such that $P$ does not have more crossings with the edges outside of $P_1\cup P_2$ than the other choice.
Since the number of crossings in the subdrawing of $P_1\cup P_2$ is trivially at most $(2d+4)^2/2\le8d^2$, the total number of crossings on the edges of $P$ is at most $\crn(G-f)/2+8d^2=:x_P$.

In case (i), draw \(f\) alongside the chosen path $P$, as in Figure~\ref{fig:slide}(a).
The newly drawn edge $f$ crosses the edges which cross $P$ also at most $x_P$ times. 
At each internal vertex of the path $P$ the drawing of $f$ can pass from one path edge to the next by crossing at most \(d-2\) incident edges in a small disc.  
After shortening possible loops along $f$, the total cost is at most \(x_P+(d-2)(d+1)<\crn(G-f)/2+9d^2\). 
This is actually the two-path argument in the embedding method of Balogh, Lea\~nos and Salazar~\cite[Section~5]{BLS}.

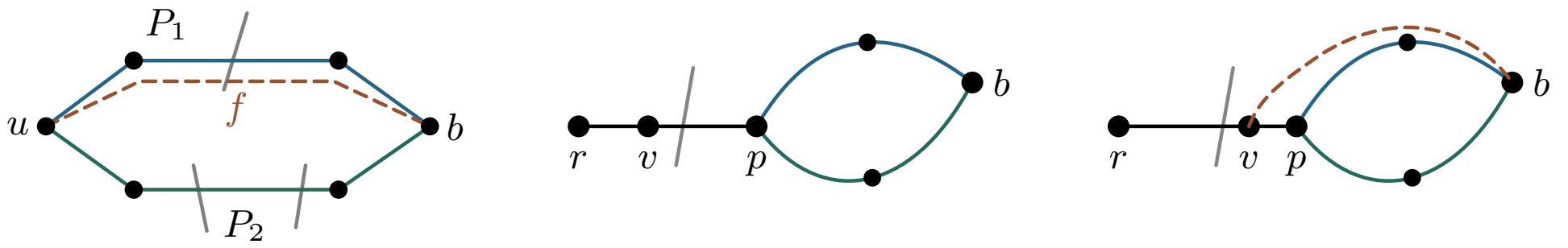
\begin{figure}[!htbp]
\centering
\begin{tikzpicture}[x=1cm,y=1cm,lab/.append style={font=\footnotesize}]
  \coordinate (u0) at (0,0); \coordinate (b0) at (3.15,0);
  \draw[edge,black!50] (1.46,.30)--(1.65,.93) (1.32,-.86)--(1.21,-.32)
                       (2.05,-.83)--(2.13,-.32);
  \draw[edge,routeblue] (u0)--(.72,.54)--(2.40,.54)--(b0);
  \draw[edge,partgreen] (u0)--(.72,-.52)--(2.40,-.52)--(b0);
  \draw[edge,warm,densely dashed] (u0)--(.77,.37)--(2.36,.37)--(b0);
  \foreach \xx/\yy in {0/0,3.15/0,.72/.54,2.4/.54,.72/-.52,2.4/-.52}{
    \node[vertex,inner sep=1.5pt] at (\xx,\yy) {};}
  \node[lab,left] at (u0) {$u$}; \node[lab,right] at (b0) {$b$};
  \node[lab,above=1pt] at (.99,.54) {$P_1$}; \node[lab,below=1pt] at (1.63,-.52) {$P_2$};
  \node[lab,text=warm] at (1.56,.14) {$f$};
  \node[panel] at (1.55,-1.2) {(a) tracing \(P_1\)};
  \foreach \shift/\vx/\id/\ttl in {4.37/.57/b/before the slide,8.80/1.07/c/after the slide}{
   \begin{scope}[xshift=\shift cm]
    \coordinate (r) at (0,0); \coordinate (v) at (\vx,0);
    \coordinate (p) at (1.46,0); \coordinate (b) at (3.23,.36);
    \draw[edge,black!48] (.80,-.32)--(.94,.48);
    \draw[edge] (r)--(p);
    \draw[edge,routeblue] (p) .. controls (1.96,.85) and (2.60,.85) .. (b) node[pos=.54,vertex,inner sep=1.45pt] {};
    \draw[edge,partgreen] (p) .. controls (1.95,-.67) and (2.71,-.62) .. (b) node[pos=.54,vertex,inner sep=1.45pt] {};
    \foreach \vtx in {r,v,p,b}{\node[vertex] at (\vtx) {};}
    \node[lab,below=2pt] at (r) {$r$}; \node[lab,below=2pt] at (v) {$v$};
    \node[lab,below=2pt] at (p) {$p$}; \node[lab,right=1pt] at (b) {$b$};
    \node[panel] at (1.58,-1.2) {(\id) \ttl};
   \end{scope}}
  \begin{scope}[xshift=8.80cm]
    \draw[edge,warm,densely dashed] (1.07,0)
       .. controls (1.16,.26) and (1.43,.43) .. (1.66,.59)
       .. controls (2.28,.98) and (2.89,.83) .. (3.23,.36);
  \end{scope}
\end{tikzpicture}
\caption{The dashed curve shows the reinserted edge. In (a) it follows the cheaper \(u\)--\(b\)
path.  In (c) it follows the cheaper  \(p\)--\(b\) path. From (b) to (c), sliding \(v\) past the grey segment leaves \(vp\) free of crossings. }
\label{fig:slide}
\end{figure}

In case (ii), see Figure~\ref{fig:slide}(b) and (c), deleting \(vb\) makes \(v\) a subdivision point on the simple arc \(rv\cup vp\).
Slide \(v\) along this arc until it is so close to \(p\) so that the short segment \(vp\) has no crossings. 
The arc composed of the union of both adjacent edges to $v$ stays the same. 
The edge $f$ is now drawn from the new position of \(v\) alongside the path $P+vp$ of length at most $d+3$, and the same argument as in (i) yields that the total cost is at most \(x_P+(d-2)(d+2)<\crn(G-f)/2+9d^2\).

So, both constructions add at most \(\crn(G-f)/2+9d^2\) crossings to the original count $\crn(G-f)$,
proving the inequality when $a_d=9d^2$. 
\end{proof}

We note in passing that \Cref{lem:traces} can be generalised in part (i) -- if there are $\gamma\geq2$ edge-disjoint $u$--$b$ paths in $G-f$. The conclusion would then read $\crn(G)\le \tfrac{\gamma+1}{\gamma}\crn(G-f)+a_{d,\gamma}$ where the additive constant now depends also on~$\gamma$.
It is not clear what should be an analogous generalisation in part (ii).

\subsection{The boundary of a small connected set}

For a vertex set \(A\), let \(N_G(A)\) be the set of its neighbours outside \(A\), and let \(\partial_G(A)\) be the set of edges with exactly one end in \(A\). 
Thus parallel edges to the same outside vertex count separately in \(\partial_G(A)\).

We shall use a simple drawing construction. A graph with \(t\) edges can be drawn in a disc with \(O(t^2)\) crossings, even if the vertices on the boundary face are prescribed. 
Place the other vertices in the interior in general position and draw straight line edges, except for parallel edge bundles where we draw them with one bend in an $\varepsilon$-neighbourhood of the straight line edge. Then \(O(t)\) edges have at most ${t\choose2}<t^2$ crossings. 

\Needspace{9\baselineskip}
\begin{lemma}
\label{lem:boundary}
For each integer \(d\ge2\) there is a constant \(b_d\) with the following property. Let \(k\ge1\)
be an integer. Suppose \(\delta(\widetilde{G})\ge3\), \(\crn(G-e)\le k-1\) for every edge \(e\), and
\(\crn(G)>\tfrac32 k+b_d\). If \(A\) is a nonempty connected set of at most \(d\) vertices, all of
degree at most \(d\) in \(G\), then
\begin{enumerate}
\item \(|N_G(A)|\ge3\),
\item if \(|\partial_G(A)|=3\), then \(A\) consists of a single cubic vertex.
\end{enumerate}
\end{lemma}
\begin{proof}
Let \(H\) consist of the edges incident with \(A\) and their endpoints, and set \(J=G-A\). The
graphs \(H,J\) partition the edges of \(G\) and meet precisely at \(N_G(A)\). Note that
\(t:=|E(H)|\le\sum_{v\in A}d_G(v)\le d^2\).

\step{At most one outside neighbour.}
Choose an edge \(f\) of \(H\), which exists since \(\delta(\widetilde{G})\ge3\). Consider the subdrawing of $J$ within an optimal drawing of \(G-f\). 
If \(A\) has no neighbour in $J$, draw all of \(H\), including \(f\), in a disc disjoint from \(J\). 
If \(A\)  has one neighbour in $J$, draw it in a small disk next to the vertex and connect it to the vertex. Either drawing gives \(\crn(G)\le k-1+t^2\) which for $b_d\geq t^2$ contradicts the assumptions (recall that $t\leq d^2$).

\step{Two outside neighbours.}
Suppose \(N_G(A)=\{u,v\}\). Let \(w\) be the maximum number of edge-disjoint \(u\)--\(v\) paths in
\(H\). Both \(u\) and \(v\) have a neighbour in the connected set \(A\), so \(w\ge1\). We shall
choose an edge \(f\) next, and then redraw \(H\) along a suitably chosen arc between $u$ and $v$.

\smallskip
\noindent\emph{Choosing the edge.}
Let \(w'\) be the maximum number of edge-disjoint \(u\)--\(v\) paths in \(H-f\) for some edge~$f$. We need
\begin{equation}
w'\ge1,\qquad \frac{w}{w'}\le\tfrac32 . \label{eq:cutratio}
\end{equation}
If \(w\ge3\), any edge \(f\) will do, since \(w'\ge w-1\). If \(w=1\) or \(2\), take a family of
\(w\) paths of minimum total length. We show that some edge \(f\) of \(H\) is unused. Deleting it
will give \(w'=w\).

Suppose the paths use every edge. One path cannot suffice, since a vertex of \(A\) on that path
would have at most two distinct neighbours, contradicting \(\delta(\widetilde{G})\ge3\). With two paths \(P_1,P_2\), the same observation
forces every vertex of \(A\) to lie on both. Orient both paths from \(u\) to \(v\). If \(a,b\)
occur in that order on \(P_1\) and in reverse order on \(P_2\), form the walks
$P_1[u,a]\,P_2[a,v]$ and $P_2[u,b]\,P_1[b,v]$,
where brackets denote subpaths. The walks are edge-disjoint and omit the segment between \(a\) and
\(b\) on each original path. Erasing cycles gives a shorter pair. The two paths must therefore
visit all vertices in the same order. Their edges can give each vertex of \(A\) only its
predecessor and successor as neighbours, contradicting \(\delta(\widetilde{G})\ge3\). 
This proves~\eqref{eq:cutratio}.

\smallskip
\noindent\emph{Redrawing \(H\).}
Take a simple optimal drawing of \(G-f\). Let \(y_J\) count crossings within \(J\) and \(z\)
crossings between \(J\) and \(H-f\). Then \(y_J+z\le k-1\). Of \(w'\) edge-disjoint \(u\)--\(v\)
paths in \(H-f\), one crosses \(J\) at most \(z/w'\) times. Erase \(H-f\) and keep the trace of
this path. By subdividing its trace at self-crossings and taking a simple path in the subdivision,
we obtain a simple arc \(\gamma\) from \(u\) to \(v\) crossing the edges of $J$ (but not the vertices) at most $z/w'$ times.

By the edge version of Menger's theorem, e.g.~\cite[Section~3.3]{Diestel}, \(H\) has a \(u\)--\(v\) cut $X\subseteq E(H)$ of size~\(w\).
Draw each side of the cut, that is each of the two sides of $H-X$, in a small free sector at $u$ and at $v$, respectively, so that these drawings are disjoint from $J$ except at $u$ and $v$.
Reconnecting the loose ends of each of the $w$ cut edges together closely along $\gamma$ in matching orders
(Figure~\ref{fig:strip}) yields a drawing of $G$ in which each edge of $X$ crosses \(J\) at most \(z/w'\) times. 
All further crossings of $H$ occur inside the sectors and cost less than $t^2$, as above.
This yields a drawing of \(G\), including \(f\), with $\crn(G)\le y_J+wz/w'+t^2\le\frac{w}{w'}(y_J+z)\le\tfrac32(k-1)+t^2$.
As \(t\le d^2\), choosing \(b_d\geq t^2\) excludes all cases with at most two outside neighbours. This proves (i).

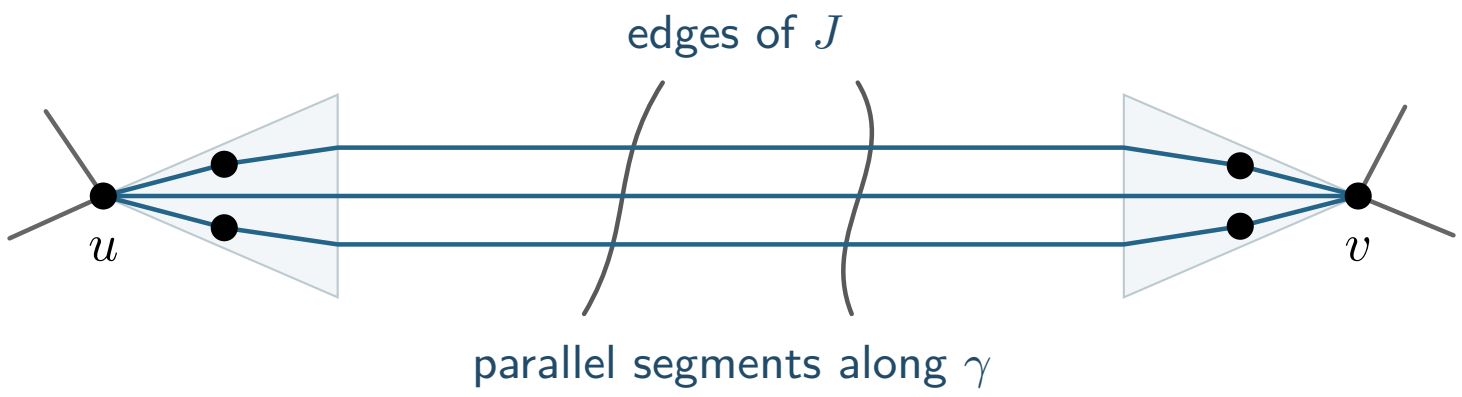
\begin{figure}[!htbp]
\centering
\begin{tikzpicture}[x=1cm,y=1cm]
 \filldraw[fill=routeblue!6,draw=lightline,thin] (0,0)--(1.55,.67)--(1.55,-.67)--cycle;
 \filldraw[fill=routeblue!6,draw=lightline,thin] (8.3,0)--(6.75,.67)--(6.75,-.67)--cycle;
 \draw[edge,black!60] (-.38,.56)--(0,0)--(-.62,-.28) (8.61,.59)--(8.3,0)--(8.93,-.26);
 \draw[edge,black!65] (3.18,-.78) .. controls (3.52,-.21) and (3.34,.19) .. (3.70,.75);
 \draw[edge,black!65] (4.95,-.78) .. controls (4.72,-.22) and (5.30,.24) .. (4.99,.75);
 \draw[edge,routeblue] (0,0)--(.8,.21)--(1.55,.32)--(6.75,.32)--(7.52,.20)--(8.3,0);
 \draw[edge,routeblue] (0,0)--(8.3,0);
 \draw[edge,routeblue] (0,0)--(.8,-.21)--(1.55,-.32)--(6.75,-.32)--(7.52,-.20)--(8.3,0);
 \foreach \xx/\yy in {0/0,.8/.21,.8/-.21,7.52/.20,7.52/-.20,8.3/0}{\node[vertex] at (\xx,\yy) {};}
 \node[lab,below=4pt] at (0,0) {$u$}; \node[lab,below=4pt] at (8.3,0) {$v$};
 \node[panel,above=1pt] at (4.18,.75) {edges of \(J\)};
 \node[panel,below=2pt] at (4.16,-.78) {parallel segments along \(\gamma\)};
\end{tikzpicture}
\caption{Each shaded sector contains one side of the cut. Here \(w=3\). Every crossing of the
guide arc with an edge of \(J\) gives one crossing on each of the three segments.}
\label{fig:strip}
\end{figure}

\step{Three boundary edges.}
Suppose \(\partial_G(A)=\{f_1,f_2,f_3\}\) and \(|A|>1\). By (i), the outside ends \(b_1,b_2,b_3\) of $f_1,f_2,f_3$ are distinct. 
The degree sum gives \(2|E(G[A])|+3\ge3|A|>2(|A|-1)+3\). Thus \(G[A]\) contains a spanning tree $T\subseteq G[A]$ and an edge \(f\in E(G[A])\setminus E(T)\).

In a simple optimal drawing of \(G-f\), consider \(T\) together with $f_1,f_2,f_3$. 
Its minimal subtree $T_1\subseteq H-f$ joining \(b_1,b_2,b_3\) consists of three branches that meet at one vertex \(a\in A\). 
Retain \(J\) and redraw all of $G[A]$, including \(f\), in a small disc around the former position of \(a\) disjoint from $J$, and reconnect the edges $f_1,f_2,f_3$ with their ends $b_1,b_2,b_3$ in $J$ along the three branches of $T_1$ in the original drawing of $H-f$.
Crossings with \(J\) along the edges $f_1,f_2,f_3$ clearly cost at most the old crossings between \(H-f\) and \(J\). 
Since within $H$ there are altogether at most $t^2$ crossings, we obtain \(\crn(G)\le \crn(G-f)+t^2\leq k-1+t^2\).
Again, for \(b_d\geq t^2\) we get a contradiction.
Hence \(|A|=1\), and its vertex has degree three because there are no loops.
\end{proof}

\section{Counting the edges to the separator}\label{sec:trees}

Fix \(\varepsilon=1/1000\) and take \(r_\varepsilon\) from Lemma~\ref{lem:separator}.
Apply Lemmas~\ref{lem:traces} and~\ref{lem:boundary} with \(d=r_\varepsilon\), and set
\(c=\max\{0,a_{r_\varepsilon},b_{r_\varepsilon}\}\). These choices
precede the choice of \(k\) and \(G\). Suppose Theorem~\ref{thm:main} fails for loopless
multigraphs with this \(c\). For an integer \(k\) for which it fails, choose \(G\) as in
Lemma~\ref{lem:min}. Lemma~\ref{lem:separator} supplies \(B\) and \(U=V(G)\setminus B\). Write
$\eta=\varepsilon n$ and $m_B=|E_G(U,B)|$,
where \(E_H(X,Y)\) denotes the edges of \(H\) between disjoint sets \(X,Y\). Thus \(|B|\le\eta\),
and \(m_B\) counts each edge copy separately. Our aim in this section is to prove \(m_B\le62\eta\).

Every nonempty connected set in \(U\) satisfies Lemma~\ref{lem:boundary}, since it lies in one
component of \(G[U]\). Both configurations of Lemma~\ref{lem:traces} are impossible in \(G\),
since either would give \(\crn(G)\le\tfrac32(k-1)+a_{r_\varepsilon}\).

Let \(h\) be the number of components of \(G[U]\). Let \(\iota_0\) count the pairs \((K,b)\) for which
the component \(K\) has a neighbour \(b\in B\). Each \(K\) has at least three such neighbours by
Lemma~\ref{lem:boundary}(i). On the other hand, applying~\eqref{eq:incidence} to the \(h\)
components gives \(\iota_0\le2h+\eta\). Therefore
$3h\le \iota_0\le2h+\eta, h\le\eta$, and $\iota_0\le3\eta$.             

\subsection{Reading the edges to \texorpdfstring{\(B\)}{B} on a tree}\label{sec:paths}

Fix a component \(K\) of \(G[U]\). Delete its bridges and call the resulting connected components
the \emph{parts} of \(K\). Each part is bridgeless, since every nonbridge lies on a cycle and a
cycle uses no bridge. A part may be a single vertex. Contracting the parts and retaining the
bridges gives a tree \(T_K\). Indeed, any cycle after contraction would give an alternative route
around one of its bridges. We refer to the vertices of \(T_K\) as \emph{nodes}.

For each edge \(ub\) with \(u\in V(K)\) and \(b\in B\), put a \emph{mark} labelled \(b\) on the
node whose part contains \(u\). A mark represents one specific edge. Two parallel edges therefore
give two occurrences of the same label. Over all trees, the number of marks is \(m_B\).

\begin{lemma}
\label{lem:nothree}
A simple path of \(T_K\), including a path consisting of one node, contains at most two
occurrences of any one label.
\end{lemma}
\Needspace{4\baselineskip}
\begin{proof}
We use one property of a bridgeless part. From a vertex \(u\) in the part one can reach any two
prescribed vertices by edge-disjoint paths, one to each. The prescribed vertices may coincide, and
either path may have length zero. To prove this, join the prescribed vertices to a new sink by one
edge each, using parallel edges if necessary. A cut separating \(u\) from the sink either splits
the part, in which case it contains at least two edges of the part, or contains the two new edges.
Menger's theorem gives the required paths.

Suppose now that three marks labelled \(b\) lie on one tree path. Choose the middle mark, breaking
ties arbitrarily. Let \(ub\) be its edge and \(A\) its part. We will reach \(b\) from \(u\) by
each of the other two marked edges (Figure~\ref{fig:three}).

\begin{figure}[!htbp]
\centering
\begin{tikzpicture}[x=1cm,y=1cm,lab/.append style={font=\small}]
  \coordinate (b) at (3.08,2.05);
  \foreach \prefix/\xx in {l/.68,m/3.08,r/5.48}{
    \coordinate (\prefix T) at (\xx,.80);
    \coordinate (\prefix L) at (\xx-.62,-.14);
    \coordinate (\prefix R) at (\xx+.62,-.14);
    \filldraw[edge,fill=routeblue!5,draw=black!48] (\prefix T)--(\prefix L)--(\prefix R)--cycle;}
  \draw[edge,routeblue] (mT)--(mL)--(lR)--(lT)--(b);
  \draw[edge,partgreen] (mT)--(mR)--(rL)--(rT)--(b);
  \draw[edge,warm,densely dashed] (mT)--(b);
  \foreach \v in {lT,lL,lR,mT,mL,mR,rT,rL,rR,b}{\node[vertex] at (\v) {};}
  \node[lab,above=1pt] at (b) {$b$}; \node[lab,below=2pt] at (mT) {$u$};
  \node[lab,below=2pt] at (mL) {$x$}; \node[lab,below=2pt] at (mR) {$y$};
  \draw[-{Stealth[length=2mm]},inkblue] (6.50,.34)--(7.65,.34);
  \node[panel,above=2pt] at (7.06,.34) {contract};
  \draw[edge] (8.13,.34)--(9.61,.34)--(11.09,.34);
  \foreach \xx in {8.13,9.61,11.09}{
    \node[block] at (\xx,.34) {}; \node[lab,above=4pt] at (\xx,.34) {$b$};}
  \node[panel] at (9.61,-.22) {three nodes of \(T_K\)};
\end{tikzpicture}
\caption{Three occurrences of \(b\) on a tree path would give two alternatives to the middle edge
\(ub\) (dashed). In this example they leave the middle triangle at \(x\) and \(y\), cross
different bridges, and finish along the two outer edges to \(b\).}
\label{fig:three}
\end{figure}
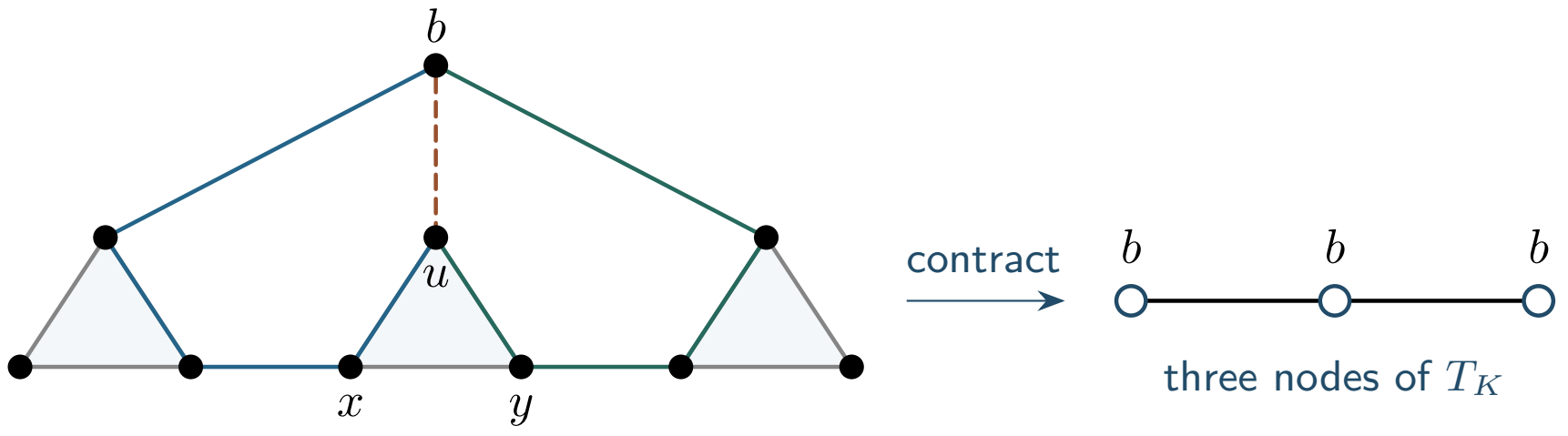

For each outer mark, choose an exit from \(A\). If that mark lies on \(A\), take the endpoint in
\(A\) of its own marked edge. Otherwise take the endpoint in \(A\) of the first bridge towards its
part. Use the property above to reach the two exits by edge-disjoint paths. Beyond \(A\), follow
the bridges of the tree path, joining successive bridge ends within the intervening parts. If both
routes leave \(A\), they go in opposite directions and use disjoint parts. Finish each route with
its own edge to \(b\). This works even if all three marks lie on \(A\).

The resulting routes avoid \(ub\). Erasing cycles gives two edge-disjoint \(u\)--\(b\) paths. All
their vertices except \(b\) lie in \(K\), so their lengths are at most \(|V(K)|\le r_\varepsilon\), and their
internal degrees in \(G\) are at most \(r_\varepsilon\). Lemma~\ref{lem:traces}(i) gives a contradiction.
\end{proof}

\subsection{One bound for both counts}\label{sec:pathcount}

Root a tree with at least two nodes at a leaf. At each node with children, keep one edge to a
child. The kept edges partition the nodes into downward paths, allowing paths with one node. Each
path ends at a different childless node, and each childless node ends a path. Hence a tree with
\(\ell_T\) leaves gives \(\ell_T-1\) paths, since the root has a child. A singleton tree gives one
path and has no leaves. Throughout, a leaf means a node of degree one.

We use this partition twice. The following lemma gives the conditions we need to check each time.

\begin{lemma}
\label{lem:pathcount}
For each \(K\), take at most one tree and partition it as above. Mark its nodes with some of the
edges from \(K\) to \(B\), using each edge once and its end in \(B\) as its label. Suppose that
\begin{enumerate}
\item each simple tree path contains at most two occurrences of each label,
\item each partition path \(P\) has a nonempty connected set \(X_P\subseteq V(K)\) containing the
\(K\)-ends of its marked edges, and all the sets \(X_P\) are pairwise disjoint.
\end{enumerate}
If there are \(q\) partition paths in total, then
$\#\{\text{marks}\}\le 2q+\eta+\iota_0\le 2q+4\eta$.
\end{lemma}
\begin{proof}
\emph{First occurrences.} Count each label \(b\) once on every partition path \(P\) on which it
occurs. The edge giving the mark joins \(X_P\) to \(b\). By~\eqref{eq:incidence}, the first
occurrences contribute at most \(2q+\eta\).

\smallskip
\Needspace{7\baselineskip}
\noindent\emph{Second occurrences.} Fix a component \(K\) and a label \(b\). Suppose two of its
partition paths each contain two marks labelled \(b\). On path \(i\), let \(x_i,y_i\) be the nodes
of the upper and lower marks, in the order away from the root. If both marks occupy one node, then
\(x_i=y_i\).

Let \(z\) be the last common node on the paths from the root to \(y_1\) and \(y_2\). The two
segments from \(x_i\) to \(y_i\) are disjoint, so at least one misses \(z\) and lies entirely
below it. The path from \(y_1\) to \(y_2\) contains that entire segment and the lower mark on the
other path (Figure~\ref{fig:paths}). These are three occurrences of \(b\), contrary to (i).

\begin{figure}[!htbp]
\centering
\begin{tikzpicture}[x=1cm,y=1cm,lab/.append style={font=\small}]
  \coordinate (r) at (0,1.70); \coordinate (x1) at (0,1.02); \coordinate (z) at (0,.29);
  \coordinate (y1) at (-1.72,-.70); \coordinate (x2) at (1.34,-.48); \coordinate (y2) at (2.76,-1.29);
  \draw[line width=5.6pt,routeblue!28,line cap=round,line join=round] (y1)--(z)--(x2)--(y2);
  \draw[edge] (r)--(x1)--(z)--(y1) (x2)--(y2);
  \draw[edge,densely dashed] (z)--(x2);
  \foreach \v in {x1,y1,x2,y2}{\node[vertex,fill=inkblue,inner sep=2.2pt] at (\v) {};}
  \foreach \v in {r,z}{\node[block,inner sep=2pt] at (\v) {};}
  \node[lab,right] at (r) {root};
  \node[lab,left] at (x1) {$x_1$}; \node[lab,left] at (z) {$z$};
  \node[lab,left] at (y1) {$y_1$}; \node[lab,above right=1pt] at (x2) {$x_2$};
  \node[lab,right] at (y2) {$y_2$};
\end{tikzpicture}
\caption{Suppose the four filled nodes carry the same label. The solid edges form two downward
paths. The path from \(y_1\) to \(y_2\) (shaded) passes through both marks \(x_2,y_2\) below the
common ancestor \(z\), as well as \(y_1\).}
\label{fig:paths}
\end{figure}
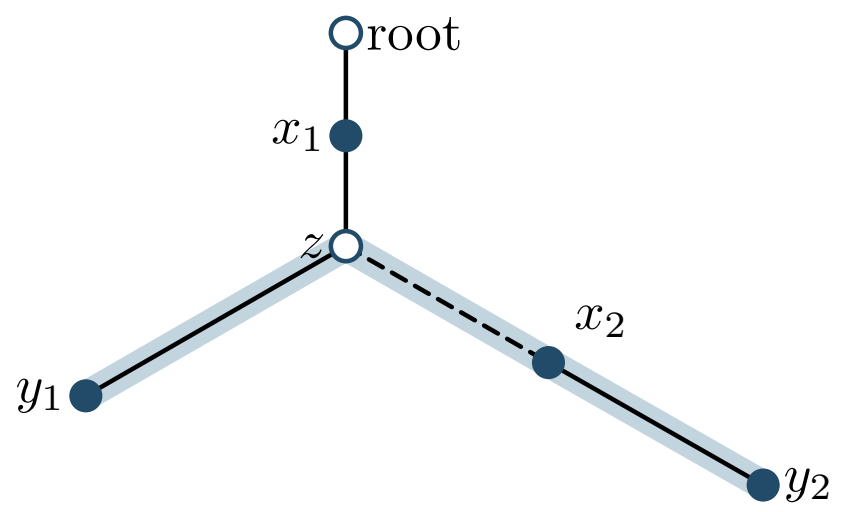

Thus a pair \((K,b)\) contributes at most one second occurrence in the entire partition. There are
at most \(\iota_0\) such pairs and no third occurrences. Adding the two counts and using
\(\iota_0\le3\eta\) proves the lemma.
\end{proof}

\subsection{The first count}\label{sec:first}

Let \(\ell\) be the total number of leaves of the trees \(T_K\). A leaf part is joined to the rest
of \(K\) by just one bridge. It has at least three distinct outside neighbours by
Lemma~\ref{lem:boundary}(i), so at least two of them belong to \(B\). Consequently \(m_B\ge2\ell\).

For the upper bound, partition the \(T_K\) into downward paths. There are at most \(\ell+h\)
paths, including one for each singleton tree. For \(X_P\), take all vertices in the parts on
\(P\). The bridges on \(P\) make \(X_P\) connected, and different paths use disjoint parts.
Lemma~\ref{lem:nothree} supplies the restriction on labels. Lemma~\ref{lem:pathcount} now gives
\(m_B\le2(\ell+h)+4\eta\). Since \(h\le\eta\), we have
$2\ell\le m_B\le2\ell+6\eta$, so $\sigma:=m_B-2\ell\le6\eta$.

Call a leaf \emph{regular} if it carries exactly two marks. Its part has three boundary edges,
namely the bridge and the two edges to \(B\). By Lemma~\ref{lem:boundary}(ii) the part is a single
cubic vertex. Its two neighbours in \(B\) are distinct.

Set aside two marks at every leaf. Exactly \(\sigma\) marks remain, including all marks on
nonleaves. An irregular leaf still has at least one mark. A node of tree degree two has only two
bridges leaving its part. Lemma~\ref{lem:boundary}(i) forces a further neighbour in \(B\), so this
node also carries a mark. It follows that
\begin{equation}
      \#\{\text{irregular leaves}\}\le\sigma,\qquad
      \#\{\text{nodes of tree degree two}\}\le\sigma .            \label{eq:accounts}
\end{equation}
The first estimate controls how many leaf pairs we can transfer. The second helps bound the number
of leaves left after pruning.

\Needspace{12\baselineskip}
\subsection{The second count}\label{sec:second}

Discard all trees with at most two nodes. 
For trees with at least three nodes, remove all its leaves at once to obtain a nonempty tree \(T'_K\). Move the two marks of each removed regular leaf to its neighbour. 
Let \(m_B'\) denote the number of transferred marks.
Each mark still represents its original edge of \(G\). 

\step{The path bound after pruning.}
We check the two hypotheses of Lemma~\ref{lem:pathcount}. First suppose that three marks labelled
\(b\) lie on a path of \(T'_K\). They came from different regular leaves, since each such leaf had
two distinct neighbours in \(B\). Let \(v\) be the cubic vertex that supplied the middle mark,
\(vp\) its bridge and \(A\) the surviving part containing \(p\). The mark has moved from \(v\) to
the node for \(A\). We shall construct the two routes from \(p\).

For each outer mark, choose an exit from \(A\) as follows. If that mark is now on \(A\), use the
end in \(A\) of its supplying leaf's bridge. Otherwise use the end in \(A\) of the first bridge
towards the mark. As in Lemma~\ref{lem:nothree}, join \(p\) to the two exits by edge-disjoint
paths. Continue through the surviving parts when necessary. If both routes continue, they go in
opposite directions on the tree path. Each route finishes along its own supplying leaf's bridge
and original edge to \(b\). Neither uses \(v\), even if all three marks now lie on \(A\). After
erasing cycles, the two \(p\)--\(b\) paths in \(G-v\) have length at most \(|V(K)|\le r_\varepsilon\). Their
internal vertices and \(p\) have degree at most \(r_\varepsilon\) in \(G\). Applying the sliding case of
Lemma~\ref{lem:traces}(ii) to \(vb\) gives the required contradiction.

For the second hypothesis, partition each \(T'_K\) into downward paths. For each path take its
surviving parts and adjoin the original regular leaves that supplied its marks, with their
bridges. The resulting vertex sets are connected and contain the original ends of their marked
edges. They are disjoint because each supplying leaf has just one surviving neighbour. If \(q'\)
is the number of partition paths after pruning, Lemma~\ref{lem:pathcount} now gives
\(m_B'\le2q'+4\eta\).

\step{How many marks survive?}
At most \(2h\le2\eta\) original leaves belonged to discarded trees, and at most \(\sigma\le6\eta\)
were irregular. Each of the remaining leaves contributes two marks. Hence
$m_B'\ge2(\ell-2h-\sigma)\ge2\ell-16\eta$.   

\step{How many paths remain?}
Write \(\ell'\) for the total number of leaves of the pruned trees. Each such leaf had exactly one
neighbour that survived. If its original degree was at least three, it lost at least two leaf
neighbours. If its original degree was two, it lost one. There are at most \(\sigma\) nodes of the
latter kind by~\eqref{eq:accounts}. No deleted leaf is counted at two surviving nodes. Thus
\(2\ell'-\sigma\le\ell\). Including the singleton pruned trees gives
$q'\le\ell'+h\le\tfrac12(\ell+\sigma)+h
\le\tfrac12\ell+4\eta$.

Substituting this bound for \(q'\) gives \(m_B'\le\ell+12\eta\). The two counts are summarised below.
\begin{center}
\small
\begin{tabular}{@{}lccc@{}}
\toprule
 & Partition paths, at most & Marks, at least & Marks, at most\\
\midrule
Original trees & \(\ell+\eta\) & \(2\ell\) & \(2\ell+6\eta\)\\
Pruned trees   & \(\ell/2+4\eta\) & \(2\ell-16\eta\) & \(\ell+12\eta\)\\
\bottomrule
\end{tabular}
\end{center}
The last two entries in the second row bound the same number \(m_B'\). Hence
\(2\ell-16\eta\le\ell+12\eta\), giving \(\ell\le28\eta\). Returning to the original trees, we obtain
\begin{equation}
m_B=2\ell+\sigma\le56\eta+6\eta=62\eta .    \label{eq:upper}
\end{equation}

\Needspace{11\baselineskip}
\section{The lower bound from chordless graphs}\label{sec:deficit}

A \emph{chord} of a cycle is an edge between two nonconsecutive vertices of the cycle. A graph is
\emph{chordless} if none of its cycles has a chord. For a (simple) graph \(H\), let
$\Def(H)=\sum_{v\in V(H)} \max\{3-d_H(v),0\}$.

We will use the below theorem by Dirac, as stated by Oxley. 
A graph is \emph{minimally $2$-connected} if it is $2$-connected and removing any single edge makes the graph $1$-connected. 

\begin{theorem}[{Dirac~\cite{Dirac}, Oxley~\cite{Oxley}}]
\label{thm:dirac}
If \(H\) is a minimally \(2\)-connected graph on \(s\) vertices, then $H$ has at least \(\lceil (s+4)/3\rceil\) vertices of  degree $2$.
\end{theorem}

\begin{lemma}
\label{lem:deficit}
Every simple chordless graph \(H\) satisfies \(\Def(H)\ge |V(H)|/5\).
\end{lemma}
\begin{proof}
We first show that a simple, connected, bridgeless, chordless graph \(H\) of order \(s\ge2\) has at least
\(2+s/5\) vertices of degree $2$. 
Since $H$ is simple, connected, and bridgeless, we know \(s\ge3\).

Suppose that \(H\) is 2-connected. 
We claim $H$ is minimally $2$-connected. 
Indeed, if \(H-uv\) is still 2-connected for an edge  $uv$, then $H-uv$ contains two internally vertex-disjoint \(u\)--\(v\)
paths. 
Since \(H\) is simple, both paths have length at least $2$, so their union is a cycle with chord \(uv\), which is a contradiction. 
Thus \(H\) is minimally 2-connected. 
By \cref{thm:dirac}, \(H\) has at least
\(\lceil(s+4)/3\rceil\) vertices of degree $2$. 
For \(s\ge3\), we have \(\lceil(s+4)/3\rceil\ge2+s/5\), so $\Def(H)\geq |V(H)|/5$. 

Now suppose $H$ is not $2$-connected, so \(H\) has a cut vertex \(v\).
Let $H_1, H_2$ be maximal connected subgraphs of $H$ where $H_1\cap H_2=\{v\}$ and $H_1\cup H_2=H$.
Note that each $H_i$ is bridgeless and has no chorded cycles. 
By induction, each $H_i$ has $2+|V(H_i)|/5$ vertices of degree $2$. 
Since $v$ cannot be a vertex of degree $2$ in $H$, we know $H$ has at least $(2+|V(H_1)|/5)+(2+|V(H_2)|/5)-2=2+(s+1)/5>2+s/5$ vertices of degree $2$. 
This proves the claim.

Now let \(H\) be any simple, chordless graph. Let $r$ be the number of bridgeless components when removing all bridges from $H$. 
Each such component \(H\) with at least two vertices contributes at least \(2+|V(H)|/5\) to $\Def(H)$. 
A single vertex component contributes \(3\ge2+1/5\). Restoring a bridge reduces the
deficit by at most $1$ at each endpoint. 
There are \(r-1\) bridges, so
$\Def(H)\ge 2r+|V(H)|/5-2(r-1)> |V(H)|/5$.
\end{proof}

\Needspace{14\baselineskip}
\begin{proof}[Proof of Theorem~\ref{thm:main}]
The simple graph \(\widetilde{G}[U]\) is chordless. 
Otherwise, choose an edge $uv$ of \(G\) that is a chord of a cycle $C$. 
The two $u-v$ paths on $C$ give the two short paths in Lemma~\ref{lem:traces}(i). Note that these two paths have bounded lengths and bounded internal degrees since $C$ is a subgraph of \(G[U]\).
This contradicts the choice of \(G\).

Each \(u\in U\) has at least three neighbours in \(\widetilde{G}\), so $u$ has at least \(\max\{3-d_{\widetilde{G}[U]}(u),0\}\) neighbours in \(B\). 
The sum of the numbers of neighbours in \(B\), over all \(u\in
U\), is exactly \(|E_{\widetilde{G}}(U,B)|\). Lemma~\ref{lem:deficit} and~\eqref{eq:upper} give
$\frac{n-|B|}{5}\le\Def(\widetilde{G}[U])\le |E_{\widetilde{G}}(U,B)|\le m_B\le62\eta$.
Since \(|B|\le\eta=n/1000\), the left side is at least \(999n/5000\), whereas
\(62\eta=310n/5000\). 
This contradiction proves
Theorem~\ref{thm:main}.
\end{proof}

\section{Conclusions}
We have proved that every $k$-crossing-critical graph, for any $k\geq1$, has crossing number at most \(\tfrac32 k+c\) for an absolute constant \(c\) independent of $k$.
Our result breaks a barrier of $2k+o(k)$ upper estimates at which research has seemed to struggle in quite many years.
Furthermore, it appears that our arguments could be generalized to provide an upper bound of order $\frac43k+o(k)$, but some parts, namely in \Cref{sec:trees}, then become very technical.
On the other hand, the central question of validity of an upper bound matching the order of $k+O(\sqrt k)$ of existing lower bounds remains wide open.

\subsection*{Declaration of LLM use}
The authors used OpenAI's models to assist with mathematical exploration and argument refinement. 
The authors themselves verified all mathematical arguments, revised the exposition, and take full responsibility for the contents of the article.

\subsection*{Acknowledgements}

This work is supported by the Center for Advanced Computation at Korea Institute for Advanced Study.

Ilkyoo Choi was supported by the Hankuk University of Foreign Studies Research Fund, the National Research Foundation of Korea (NRF) grant funded by the Korea government (MSIT) (RS-2025-23324220), Institute for Basic Science (IBS-R029-C1), and the Korea Institute for Advanced Study (KIAS) grant funded by the Korea government.

Tomáš Masařík was supported by the Polish National Science Centre SONATA-
17 grant number 2021/43/D/ST6/03312.

We would like to thank the organisers of the Crossing Numbers
Workshop 2026 where this research was initiated.

\bibliography{references}
\bibliographystyle{alphaurl}

\iffalse

\fi
\end{document}